\documentclass[12pt]{amsproc}
\usepackage{amssymb,latexsym}
\theoremstyle{plain}
\newtheorem{theorem}{Theorem}[section]
\newtheorem{corollary}[theorem]{Corollary}
\newtheorem{example}[theorem]{Example}

\newtheorem{proposition}[theorem]{Proposition}

\theoremstyle{definition}
\newtheorem{definition}[theorem]{Definition}

\theoremstyle{remark}
\newtheorem{remark}[theorem]{Remark}

\begin{document}

\title{Projective multiresolution analyses for dilations in higher dimensions}
\author[J.A.~Packer]{Judith A. Packer}
\address{Department of Mathematics, University of Colorado, CB 395, Boulder, Colorado, 80309-0395, U.S.A.}
\email{packer@euclid.colorado.edu}
\keywords{module frames, finitely generated projective modules, dilations,
$K$-theory, wavelets, C*-algebras, Hilbert C*-module}
\subjclass{Primary 46L99; Secondary 42C15, 46H25, 19M05}
%%%%%%%%%%%%%%%%%%%%%%%%%%%%%%%%%%%%%%%%%%%%%%%%%%%%%%%%%%%%%%%%%%%%%%%%%%%%%
\begin{abstract}
We continue the study of projective module wavelet frames corresponding to diagonal dilation matrices on $\mathbb R^n$ with integer entries, focusing on the construction of a projective multi-resolution analysis corresponding to dilations whose domains are finitely generated projective modules over continuous complex-valued functions on $n$-tori, $n\geq 3$. We are able to generalize some of these results to dilation matrices that are conjugates of integral diagonal dilation matrices by elements of $SL(n,\mathbb Z).$ We follow the method proposed by the author and M. Rieffel, and are able to come up with examples of non-free projective module wavelet frames which can be described via this construction. As an application of our results, in the case $n=3,$ when the dilation matrix is a constant multiple of the identity, we embed every finitely generated module as an initial module. 
\end{abstract}
%%%%%%%%%%%%%%%%%%%%%%%%%%%%%%%%%%%%%%%%%%%%%%%%%%%%%%%%%%%%%%%%%%%%%%%%%%%%%%%
\maketitle 

%%%%%%%%%%%%%%%%%%%%%%%%%%%%%%%%%%%%%%%%%%%%%%%%%%%%%%%%%%%%%%%%%%%%%%%%%%%%%%%
\section{Introduction}

Let $A$ be an $n\times n$ integer dilation matrix.  In a previous paper \cite{pacrief2}, the author and M. Rieffel studied the existence of projective multiresolution analyses corresponding to dilation by certain $2 \times 2$ matrices $A$ in a Hilbert $C(\mathbb T^n)$-module $\Xi,$  where we recall that for $n\in\;\mathbb N,\;n\geq 2,$  the right $C(\mathbb T^n)$ module $\Xi$ is defined as the completion of $C_c(\mathbb R^n)$ under the norm determined by the following $C(\mathbb T^n)$-valued inner product:
\begin{equation}
\label{eq ip2}
    \langle \xi, \eta \rangle_{C(\mathbb T^n)}(t) := \sum_{p \in {\mathbb Z}^n}
                   (\overline{\xi}\eta)(t-p).
\end{equation}
for $t \in {\mathbb R}^n.$ 
We remark that space $\Xi$ was first introduced by G. Zimmermann in Chapter V of his thesis \cite{Z1},  where the notation $L^{2,\infty}(\mathbb Z^n,\mathbb T^n)$ was used, and was studied further by J. Benedetto and Zimmermann in \cite{BZ}. 

One verifies (c.f. \cite{pacrief2}, Section 1) that $\Xi$ consists of all bounded continuous functions $\xi$ on $\mathbb R^n$ for which there is a constant, $K$, such that $\sum_{p\in\mathbb Z^n}|\xi(x-p)|^2\leq K$ for each $x\in\mathbb R^n.$ The right module action of $C(\mathbb T^n)$ on $\Xi$ is given by pointwise multiplication. 
One also checks that 
$$\|\langle \xi, \xi \rangle_{C(\mathbb T^n)}\|\;\geq\;\int_{\mathbb R^n}|\xi(t)|^2\;dt,$$
so that $\Xi\;\subseteq\;L^2(\mathbb R^n).$ 
In \cite{pacrief2}, Proposition 5, sufficient decay conditions were provided on any continuous $\xi\;\in\;L^2(\mathbb R^n)$ for $\xi$ to be an element of $\Xi.$ 

   In a previous paper \cite{pacrief2}, the notion of  a projective multiresolution analysis for a general $n\times n$ dilation matrix was defined as follows: 
\begin{definition}
\label{def qmag}
[\cite{pacrief2}, Definition 4] Fix $n\;\in\;\mathbb N,$ let $A$ be a $n\times n$ integer dilation matrix, and let $\Xi$ be the right-rigged Hilbert $C(\mathbb T^n)$ module defined above.  A sequence $\{V_j\}_{j\in\mathbb Z}$ of subspaces of $\Xi$ is called a {\bf projective multiresolution analysis} for dilation by $A$ if:
\renewcommand{\labelenumi}{(\roman{enumi})}
\begin{enumerate}
\item $V_0$ is a finitely generated projective $C(\mathbb T^n)$ submodule of $\Xi;$
\item $V_j\;=\;D^{j}(V_0),\;\forall\;j\;\in\;\mathbb N,$ where $D$ is as defined below,
\item $V_j\;\subset\; V_{j+1},\;\forall\;j\;\in\;\mathbb Z,$
\item $\cup_{j=0}^{\infty}V_j$ is dense in $\Xi,$ in the Hilbert $C(\mathbb T^n)$-module topology,
\item $\cap_{-\infty}^{\infty}V_j\;=\;\{0\}.$
\end{enumerate}
\end{definition}
Here $D$ is defined to be the Fourier transformed version of $D_A,$
$$D\;=\;{\mathcal F}\circ D_A \circ {\mathcal F}^{\ast},$$
where the dilation operator $D_A$ is defined on $L^2(\mathbb R^n)$ by
$$D_A(\zeta)(x)\;=\;|\text{det}(A)^{1/2}|\;\zeta(A(x)),\; \zeta\in L^2(\mathbb R^n),$$ and the Fourier transform ${\mathcal F}$ is defined on $L^2(\mathbb R^n)$ by 
$${\mathcal F}(f)(x)\;\;\int_{\mathbb R^n}f(t)\overline{e}(t\cdot x)\;dt.$$
Here $e$ is the exponential function defined on $\mathbb R$ by $e(r)=e^{2\pi i r}.$
An easy calculation shows that 
$$D\;=\;D_{B^{-1}},$$
for $B=A^t.$

It turns out that condition (v) in the definitions above is implied by conditions (i) and (ii) (\cite{pacrief2}, Proposition 13).

In \cite{pacrief2}, Theorem 6, specializing to the case where $n=2,$ for any fixed diagonal dilation matrix $A$, and any choice of finitely generated projective $C(\mathbb T^2)$-module $X(q,a),$ a projective multiresolution analysis for that dilation matrix was constructed with $V_0\;\cong \;X(q,a).$
In this situation, the isomorphism class of each of the higher dimensional modules $V_j\;=\;D^j(V_0),\;j\;\geq\;1$ was provided, and under appropriate circumstances, we gave a constructive approach which gave module frames for $\Xi$ analogous to dilates and translates of $L^2$-wavelets.

It is the aim of this paper to continue developing the methods outlined in \cite{pacrief2}, applying the main ideas to construct $C(\mathbb T^n)$-module frames for $\Xi$ in dimensions $n\;\geq \; 3.$  We discuss projective multiresolution analyses corresponding to integer-valued dilation matrices which are similar to diagonal matrices via conjugation by an element in $SL(n,\mathbb Z).$ It is interesting to note that in \cite{pacrief2}, in all of the examples of projective multi-resolution analyses corresponding to $2\times 2$ dilation matrices with positive determinant, where the  ``initial module" $V_0$ was not free,  the wavelet module $W_0\;=\;V_1\ominus V_0$ was still free.  Hence the module generators of $W_0$ in this case corresponded to multi-wavelets in the usual sense.  In contrast, in this paper an example will be given of a $3\times 3$ dilation matrix with positive determinant, and a corresponding projective multiresolution analysis where $V_0,\; V_1,$ and the wavelet module $W_0$ are all finitely generated {\bf non-free} $C(\mathbb T^3)$-modules.  Also, in the $3\times  3$ case, where the dilation matrix in question is an integer multiple of the identity matrix, we are able to construct a projective multiresolution analysis for that dilation matrix whose initial module is isomorphic to ${\mathcal P},$ where here ${\mathcal P}$ represents an arbitrary finitely generated projective $C(\mathbb T^3)$ -module.

In Section $2,$ the class of allowable ``scaling functions" which lie in $\Xi$ is widened, in Section $3$ we review $K$-theory and finitely generated projective modules for $\mathbb T^n,$ and discuss the extent to which the two-dimensional case can be generalized to higher dimensions.  In Section $4,$ we discuss the construction of projective multiresolution analyses for dilation matrices that are similar via an element of $SL(n,\mathbb Z)$ to a integral diagonal dilation matrix, and in Section $5$ a method is given for constructing projective module frames for $\Xi$ from a projective multiresolution analysis whenever module frames for $V_0$ and $W_0$ are known.  In particular, this result  shows that $W_0$ is a free module, then $W_i$ is free for all $i\;\in\;\mathbb N.$ 
We finish by posing some open questions.

The author would like to thank Marc Rieffel and Dana Williams for many useful conversations and suggestions about the results in this paper.
%%%%%%%%%%%%%%%%%%%%%%%%%%%%%%%%%%%%%%%%%%%%%%%%%%%%%%%%%%%%%%%%%%%%%%%%%%%%%%%
\section{Constructions of scaling functions whose Fourier transforms are non-compactly supported elements of $\Xi$}

  In the construction of projective multiresolution analyses given in \cite{pacrief2}, continuous Meyer-type low-pass filters corresponding to dilation by an integer $d$ with $|d|\;>\;1$ were used, whose corresponding scaling functions had continuous Fourier transforms with compact support, hence were elements of $\Xi.$   The tensor product construction then gave scaling functions for diagonal dilation matrices whose Fourier transforms were continuous and compactly supported, thus in $\Xi.$ The aim of this construction was to use these filters to aid in building the projective multiresolution analyses. However in practice, one often comes across scaling functions whose Fourier transforms $\Phi$ are not compactly supported, and it is useful to have a weaker condition that guarantees they lie in $\Xi.$ It also is of interest to identify more precisely the $C(\mathbb T^n)$-module structure of $\Xi.$ We do this in the following theorem.

\begin{theorem}
\label{thm freqscal} The Hilbert $C(\mathbb T^n)$-module $\Xi$ is isomorphic to the standard $C(\mathbb T^n)$-module $l^2(C(\mathbb T^n))\;=\;\{(f_j)_{j=1}^{\infty}\;:\; f_j\in C(\mathbb T^n),\;\sum_{j=1}^{\infty}|f_j(x)|^2\;\mbox{converges in norm in}\;C(\mathbb T^n)\}.$  In addition, if $\Phi$ is a continuous complex-valued function defined on $\mathbb R^n,$  then $\Phi\in\Xi$ if and only if the sum 
$$\sum_{p\;\in\;\mathbb Z^n}|\Phi(x-p)|^2$$
converges uniformly on the cube $[0,1]^n.$ 
\end{theorem}
\begin{proof}

Let ${\mathcal H}$ be the Hilbert $C^{\ast}(\mathbb Z^n)$-module obtained by completing the $C^{\ast}(\mathbb Z^n)$-module $C_C(\mathbb R^n),$ where we give $C_C(\mathbb R^n)$ a right $C_C(\mathbb Z^n)$-action defined by 
  \[
f\;\ast\;a\;=\;\sum_{v\in\mathbb Z^n}f(t-v)a(v),\;f\in\;C_c(\mathbb R^n),\;a\;\in\;C_c(\mathbb Z^n),
  \]
and a $C^{\ast}(\mathbb Z^n)$-valued inner product defined by
  \[
\langle f_1, f_2 \rangle_{C^{\ast}(\mathbb Z^n)}(v)=\int_{\mathbb R^n}\overline{f_1(x)}f_2(x-v)dx, f_1,f_2\;\in C_c(\mathbb R^n).
  \]

The module ${\mathcal H}$ was studied in much greater generality for any closed subgroup $H$ of the locally compact group $G$ by M. Rieffel in 1974 (\cite{Rie0}), where he studied induction of group representations by algebras. In this case $H=\mathbb Z^n$ and $G=\mathbb R^n$ and ${\mathcal H}$ serves as a strong Morita equivalence bimodule between $C^{\ast}(\mathbb Z^n)$ and the transformation group $C^{\ast}$-algebra $C_0(\mathbb R^n/\mathbb Z^n)\rtimes \mathbb R^n.$ It then follows from the work of P. Green (\cite{Gr}, Theorem 2.4) that ${\mathcal H}$ is isomorphic as a right $C^{\ast}(\mathbb Z^n)$-module to the standard $C^{\ast}(\mathbb Z^n)$-module $l^2(C^{\ast}(\mathbb Z^n)).$

One easily checks that $\Xi={\mathcal F}({\mathcal H})$ and that 
$${\mathcal F}(f\;\ast\;a )\;=\;{\mathcal F}(f)\cdot{\mathcal F}(a),\; f\in {\mathcal H},\; a\in C^{\ast}(\mathbb Z^n)$$ and 
$${\mathcal F}(\langle f_1, f_2 \rangle_{C^{\ast}(\mathbb Z^n)})\;=\;\langle {\mathcal F}(f_1), {\mathcal F}(f_2) \rangle_{C(\mathbb T^n)},$$ where here, by abuse of notation, we also denote the Fourier transform from $C^{\ast}(\mathbb Z^n)$ to $C(\mathbb T^n)$ by ${\mathcal F}.$ Hence $\Xi$ is isomorphic to ${\mathcal H}.$ But we have seen that ${\mathcal H}$ is isomorphic to $l^2(C^{\ast}(\mathbb Z^n))$ as a $C^{\ast}(\mathbb Z^n)$-module.  Thus $\Xi$ is isomorphic to the standard $C(\mathbb T^n)$-module $l^2(C(\mathbb T^n))$ as a $C(\mathbb T^n)$-module. 

As for the second statement, suppose that $\Phi$ is continuous on $\mathbb R^n$ and 
$\sum_{p\;\in\;\mathbb Z^n}|\Phi(x-p)|^2$ converges uniformly the cube $[0,1]^n.$ Since each term $|\Phi(x-p)|^2$ is continuous on $[0,1]^n,$ by the theory of uniform convergence, the sum $\sum_{p\;\in\;\mathbb Z^n}|\Phi(x-p)|^2$ is continuous on the compact cube $[0,1]^n,$ hence is bounded from above by some positive constant $K.$  But then for any $x\in \mathbb R^n,$ there exists $m\in\mathbb Z^n$ with $x'=x+m\in[0,1]^n,$ so that 
$$\sum_{p\;\in\;\mathbb Z^n}|\Phi(x-p)|^2$$
$$=\;\sum_{p\;\in\;\mathbb Z^n}|\Phi(x'-m-p)|^2\;\leq\;K,$$ and
hence $\Phi\in \Xi.$

The converse follows directly from Proposition 1 of \cite{pacrief2}. 
\end{proof}

\begin{remark}
We note that the construction of the isomorphism between ${\mathcal H}$ and $l^2(C^{\ast}(\mathbb Z^n))$ from \cite{Gr} depends on the choice of a Borel section $\alpha:\;\mathbb R^n/\mathbb Z^n\;\rightarrow \mathbb R^n.$  Thus there is no canonical choice of isomorphism between $\Xi$ and $l^2(C(\mathbb T^n)).$ 
\end{remark} 

We give an example of the use of the above theorem to show that the Fourier transform of the Haar scaling function in dimension $1$ for dilation by an integer $d>1$ is an element of $\Xi.$ Some ideas behind this construction can be found in Section $4$ of the reference \cite{gopbur}.
Fix $d\;\in\;\mathbb N,\;d>1,$ and find an orthonormal basis for $\mathbb R^d,$ with the initial vector fixed as 
$\vec{r_0}=(\frac{1}{\sqrt{d}},\frac{1}{\sqrt{d}},\cdots,\frac{1}{\sqrt{d}}),$ and the other vectors as $\vec{r_1},\cdots,\vec{r_{d-1}}.$  For $0\;\leq\;l\;\leq\;d-1,$ write
$$\vec{r_l}\;=\;(a_{l,0},a_{l,1},\cdots,a_{l,d-1}).$$
Let $\vec{v(x)}$ denote the row vector whose entries consist of the following $d$-periodic functions defined on $\mathbb R:\;(1,e(\frac{x}{d}),e(\frac{2x}{d}),\cdots,e(\frac{(d-1)x}{d}).$ 
Now define 
$$\mu_l(x)\;=\; \vec{r_l}\cdot\vec{v(x)},\;0\;\leq\;l\;\leq\;q-1,$$
where the ``$\cdot$" denotes dot product, so that
$$\mu_l(x)\;=\; \sum_{j=0}^{d-1}a_{l,j}e(\frac{jx}{d}),\;0\;\leq\;l\;\leq\;d-1.$$
The functions $\mu_i(x),\;0\leq\;i\;\leq\;d-1$
are continuous and periodic modulo $d,$ and 
$$\mu_0(x)\;=\; \sum_{j=0}^{d-1}e(\frac{jx}{d}),$$
so that $\mu_0(0)=1.$ 
As in \cite{gopbur}, one easily calculates by the law of characters on $\mathbb Z/d \mathbb Z$ that:
$$\langle \mu_l, \mu_k \rangle_{C(\mathbb R/\mathbb Z)}\;=\;\sum_{n=0}^{d-1}d\;a_{l,n}a_{k,n}$$
$$\;=\;d\delta_{l,k},$$
since the vectors $\vec{r_0},\;\vec{r_1},\;\cdots,\vec{r_{d-1}}$ form an orthonormal basis for $\mathbb R^d.$ (Here the use of the inner product notation refers to the definition given in Equation \ref{eq ip2}.)

Defining $m_i(x)\;=\;\mu_i(dx),\;0\leq\;i\;\leq\;d-1,$
it is clear that the functions $\{m_i:0\leq\;i\;\leq\;d-1\}$ will be continuous and periodic modulo $1.$ They will also  satisfy the two conditions 
$$m_0(0)=\sqrt{d},$$ and 
$$\langle m_l ,m_k \rangle_{C(\mathbb R/\frac{1}{d}\mathbb Z)}(x)\;=\;d\delta_{lk}(x).$$
The functions $\{m_j:\;0\;\leq\;j\;\leq\;d-1\}$  are all trigonometric polynomials, and $m_0$ satisfies Cohen's condition as stated in \cite{Str}.
In the usual way one uses the telescoping property to prove that
when $m_0$ is the trigonometric low-pass filter defined above, the corresponding scaling function viewed in the frequency domain will take on the closed form  
\begin{equation}
\label{eqn scale}
   \Phi(x)= \prod_{j=1}^{\infty}[\frac{1}{\sqrt{d}}m_0(\frac{x}{d^j})]\;=\;\frac{e(x)-1}{2\pi ix},
\end{equation}
which is familiar from ordinary wavelet theory as being the Fourier 
transform of the Haar scaling function. By an application of Theorem \ref{thm freqscal}, we immediately see that $\Phi$ is an element of the Hilbert $C(\mathbb T)$-module $\Xi.$ 
In the same way one can show that the Fourier transform for the Haar scaling function corresponding to dilation by $-d$ for $d\;\in\;\mathbb N,\;d\;>\;1,$ is also in $\Xi.$

%%%%%%%%%%%%%%%%%%%%%%%%%%%%%%%%%%%%%%%%%%%%%%%%%%%%%%%%%%%%%%%%%%%%%%%%%%%%%%%
\section{Projective multi-resolution analyses in higher dimensions: Diagonal matrices with integer entries}

In this section, we first discuss finitely generated projective modules over $C(\mathbb T^n),\;n\geq 3,$ and then go on to construct projective multi-resolution analyses in higher dimensions. We also discuss the relationship  of these projective modules to the $K$-groups of $C(\mathbb T^n).$ The modules we consider are adapted from $C(\mathbb T^2)$-modules first  constructed in  Section 3 of \cite{Rie1}, and were also reviewed in Section 3 of \cite{pacrief2}.  For completeness, we briefly describe these $C(\mathbb T^2)$-modules here. Throughout this section, unless otherwise specified, we view functions on $\mathbb T^n$ as functions defined on $\mathbb R^n$ which are periodic modulo $\mathbb Z^n.$
\begin{theorem}
For $q\;\in\;\mathbb N$ and $a\;\in\;\mathbb Z,$ let $X(q,a)$ denote 
the right $C(\mathbb T^2)$-module consisting of the space of 
continuous complex-valued functions $f$ on $\mathbb T \times \mathbb 
R$ which satisfy
$$f(s,t-q)\;=\;e(as)h(s,t),$$
with module action given by
$$f\cdot F(s,t)\;=\;f(s,t)F(s,t),$$
for $f\;\in\;X(q,a)$ and $F\;\in\;C(\mathbb T^2).$
Then $X(q,a)$ is a finitely generated, projective $C(\mathbb 
T^2)$-module; following the results from \cite{Rie1}, we say that 
$X(q,a)$ has dimension $q$ and twist $-a.$ The set 
$\{X(q,a):\;q\;\in\;\mathbb N,\;a\;\in\;\mathbb Z\}$ parametrizes the 
isomorphism classes of finitely generated projective $C(\mathbb 
T^2)$-modules, in the sense that if $X$ is a finitely generated 
projective $C(\mathbb T^2)$-module, there exist unique values of $q$ 
and $a$ such that $X\;\cong\;X(q,a).$ Moreover, if $q_1,\;q_2\;\in\;\mathbb N$ and 
$a_1,\;a_2\;\in\;\mathbb Z$ are fixed, then 
$$X(q_1,a_1)\;\oplus\;X(q_2,a_2)\;\cong\;X(q_1+q_2,a_1+a_2)$$ as
finitely generated projective $C(\mathbb T^2)$-modules.  Furthermore, cancellation holds for all finitely generated $C(\mathbb T^2)$-modules; that is, if
$$X(q_1,a_1)\;\oplus\;X(q_2,a_2)\;\cong\;X(q_1,a_1)\;\oplus\;X(q_3,a_3),$$
then $X(q_2,a_2)\;\cong\;X(q_3,a_3),$ so that $q_2=q_3$ and $a_2=a_3.$
\end{theorem}
\begin {proof}
For the proof of this result, refer to \cite{Rie1}, Theorem 3.9.
\end{proof}

It follows that for any positive integer $q,\; X(q,0)$ is isomorphic to the free 
$C(\mathbb T^2)$-module on $q$ generators; moreover, we obtain the well-known result that $K_0(C(\mathbb 
T^2))\;\cong\;\mathbb Z^2\;=\;\{[q,a]:\;q,a\;\in\;\mathbb Z\}$ (\cite{Rie1}).

The above analysis for finitely generated projective modules over $C(\mathbb T^2)$ extends to finitely generated projective modules over $C(\mathbb T^n).$  By definition, the stable isomorphism classes of finitely generated projective modules over $C(\mathbb T^n)$ give the positive cone of the abelian group $K_0(C(\mathbb T^n))\;\cong\;K^0(\mathbb T^n),$ hence generate the entire group $K^0(\mathbb T^n),$ where $K^0(\mathbb T^n)$ represents the Grothendieck group from stable isomorphism class of complex vector bundles over $\mathbb T^n.$  It is well-known that the ring $K^\ast(\mathbb T^n)\cong K_{\ast}(C(\mathbb T^n))$ is isomorphic to the
exterior algebra over $\mathbb Z$ on $n$ generators, 
$\bigwedge^\ast_{\mbox{\scriptsize \bf Z}}\{e_1,...,e_n\}$ 
(\cite{Ell}). Here the Chern character $ch:K_{\ast}(C(\mathbb T^n)) \rightarrow \check{H}^\ast
(\mathbb T^n,\mathbb Q)$ is integral and gives an isomorphism $ch_0:K_0(C(\mathbb T^n))
\rightarrow \check{H}^{\mbox{\scriptsize even}}(\mathbb T^n,\mathbb Z)$ and
$ch_1:K_1(C(\mathbb T^n))\rightarrow \check{H}^{\mbox{\scriptsize odd}}(\mathbb T^n,\mathbb Z)$ 
(see \cite{Ji} for a proof), and $\check{H}^\ast(\mathbb T^n,\mathbb Z)$ under cup
product, is known to be isomorphic to the exterior algebra  
$\bigwedge^\ast_{\mathbb Z}\{e_1,...,e_n\}$ 
with $\check{H}^k(\mathbb T^n,\mathbb Z)\cong
\bigwedge^k_{\mathbb Z}\{e_1,...,e_n\}$. For future reference we discuss some calculations involving the Chern character on $K_{\ast}(C(\mathbb T^n))$ under the identification with the exterior algebra on $n$ generators discussed above; see the reference \cite{paclee} for more details.  

Recall from \cite{pacrief2} that the Chern character on $K_0(C(\mathbb T^2))^+$ maps $[X(q,a)]$ to $q-ae_1\wedge e_2\;\in\;\bigwedge^{\mbox{\scriptsize even}}_{\mathbb Z}
\{e_1,e_2\}.$   Also recall that if $A\;=\;\left(\begin{array}{rr}
b&c\\
d&f
\end{array}\right)\;\in\;GL(2,\mathbb Z)$ is viewed an automorphism of $\mathbb T^2$ by the formula 
$$A(z,w)\;=\;(z^bw^c,z^dw^f),\;(z,w)\;\in\;\mathbb T^2,$$
then the corresponding map in $K$-theory, $A^{\ast}\;:\;K_0(C(\mathbb T^2))\;\rightarrow\;K_0(C(\mathbb T^2))$ is given by 
\begin{equation}
\label{eq act}
A^{\ast}(q-ae_1\wedge e_2)\;=\;q-{\mbox{det}}(A)a e_1\wedge e_2.
\end{equation}

More generally, viewing $C(\mathbb T^n)$ as the $C^{\ast}$-algebra of all continuous functions on the $n$ real variables $s_1,\;s_2,\;\cdots,s_{n-1},\;t,$ and fixing $q\;\in\;\mathbb N$ and $a_1,a_2,\cdots,a_{n-1}\;\in\;\mathbb Z,$ let 
$X(q,a_1,a_2,\cdots,a_{n-1})$ be the $C(\mathbb T^n)$ module consist of all continuous functions $f:\mathbb R^n\rightarrow\mathbb C$ satisfying the identity 
$$f(s_1,s_2,\cdots,s_{n-1},t-q)=e(a_1s_1+\cdots a_{n-1}s_{n-1})f(s_1,s_2,\cdots,s_{n-1},t),$$
It can be shown using methods from section 3 of \cite{BP2} and \cite{paclee} that $X(q,a_1,a_2,\cdots,a_{n-1})$ is a finitely generated projective $C(\mathbb T^n)$-module whose first Chern class corresponds to the element $-\sum_{1\leq k\leq n-1}a_ke_k\wedge e_n\in
\bigwedge^2_{\mathbb Z}\{e_1,...,e_{n}\}\cong H^2(\mathbb T^n,
\mathbb Z).$  For general $n,$ this method does not produce all finitely generated projective modules of $C(\mathbb T^n)$ as in general, the positive cone of $K_0(C(\mathbb T^n))$ seems to be unknown; see Exercise 6.10.2 of \cite{Bla}. Also, as indicated in \cite{Bla} and demonstrated in Section 4 of \cite{pacrief}, cancellation of projective modules over $C(\mathbb T^n)$ fails for $n\geq 5,$ although cancellation holds for $n\leq\;4.$ However, the failure of cancellation to hold does not effect calculations at the level of the $K$-groups, and combined with the natural automorphisms of $\mathbb T^n$ provided by elements in $GL(n,\mathbb Z)$ one can write down specific formulas for a large class of examples of elements from $K_0(C(\mathbb T^n))^+.$ In the special cases when $n=2$ and $n=3$ we can write down all of them. We can also write down analogues of Equation (\ref{eq act}) which tell us how certain elements of $GL(n,\mathbb Z)$ affect elements in $K_0(C(\mathbb T^n))^+,$ by using the exterior algebra methods outlined in the proof of Equation (\ref{eq act}).  We will use these facts in what follows when we construct our projective multi-resolution analyses for certain higher dimensional dilation matrices.

We now discuss the construction of projective multi-resolution analyses where the initial space $V$ is a finitely generated projective module over $C(\mathbb T^n),$ for $n\;\geq\;3.$  For simplicity, we will concentrate on the case where the module $V$ in question is of the form 
$X(q,0,0,\cdots,a_j,0,\cdots,0)$ as defined above, for $q\;\in\;\mathbb N,\;1\;\leq\;j\;\leq\;n-1$ and $a_1,\;a_2,\;\cdots,\;a_{n-1}\;\in\;\mathbb Z.$  In this case $X(q,0,0,\cdots,a_j,0,\cdots,0)$ is the $C(\mathbb T^n)$ module consisting of all continuous complex-valued functions from $\mathbb R^n$ to $\mathbb C$ which are periodic modulo $\mathbb Z$ in the first $n-1$ variables and which satisfy the identity
$$f(s_1,s_2,\;\cdots,s_{n-1},\;t-kq)\;=\;e(ka_js_j)f(s_1,s_2,\;\cdots,s_{n-1},\;t),$$
$$s_1,\;s_2,\;\cdots,\;s_{n-1},\;t\;\in\;\mathbb R,\;k\;\in\;\mathbb Z,$$
and where the action of $C(\mathbb T^n)$ on $X(q,0,0,\cdots,a_j,0,\cdots,0)$ is by pointwise multiplication. Here, as previously, $j$ represents a fixed integer between $1$ and $n-1.$
\newline We first note that by our discussion in Section 2, the module $X(q,0,0,\cdots,0,a_j,0,\cdots,0)$ is exactly the pullback of the $C(\mathbb T^2)$-module $X(q,a_j)$ obtained from the projection of $\mathbb T^n$ onto $\mathbb T^2$ given by projection in the $j^{\mbox{th}}$ and $n^{\mbox{th}}$ coordinates, and therefore, upon identifying $K_0(C(\mathbb T^n))$ with $\bigwedge^{\mbox{\scriptsize even}}_{\mathbb Z}
\{e_1,...,e_n\}$ as in section two, $[X(q,0,0,\cdots,0,a_j,0,\cdots,0)]_{K_0(C(\mathbb T^n))}$ can be identified with the element $q\;-\;a_je_j\wedge\;e_n\;\in\;\bigwedge^{\mbox{\scriptsize even}}_{\mathbb Z}
\{e_1,...,e_n\}.$ The $C(\mathbb T^n)$-valued inner product on $X(q,0,0,\cdots,0,a_j,0,\cdots,0)$ is given by 
$$\langle h_1,h_2 \rangle_{C(\mathbb T^n)}\;=\;\sum_{k=0}^{q-1}{\overline {h_1(s_1,s_2,\cdots,s_{n-1},t-k)}}\;\cdot\;h_2(s_1,s_2,\cdots,s_{n-1},t-k).$$
We now remark that it is possible to generalize Theorem 4 of \cite{pacrief2} to this $n$-dimensional setting, and embed $X(q,0,0,\cdots,0,a_j,0,\cdots,0)$ as a finitely generated projective module in the Hilbert $C(\mathbb T^n)$-module $\Xi\;\subseteq\;L^2(\mathbb R^n)$ described in Section 1, as follows:
\begin{theorem}
\label{thm quasi3}
Fix $n\;\in\;\mathbb N,$ and let $A$ be the $n\times n$ diagonal matrix 
$$A\;=\;\left(\begin{array}{rrrr}
d_1&0&\cdots&0\\
0&d_2&\cdots&0\\
\cdot&\cdot&\cdots&\cdot\\
0&0&\cdots&d_n
\end{array}\right),$$ where $d_1,\;d_2,\cdots,\;d_n\;\in\;\mathbb N$ and
$|d_j|>1,\;1\;\leq\;j\;\leq\;n.$ Fix $q\;\in\;\mathbb N$ and $a_j\;\in\;\mathbb Z,$ and let $X(q,0,0,\cdots,0,a_j,0,\cdots,0)$ be the $C(\mathbb T^n)$-module defined in Section 2.
Then there exists a $C(\mathbb T^n)$-module monomorphism ${\mathcal R}:\;X(q,0,0,\cdots,0,a_j,0,\cdots,0)\;\rightarrow\;\Xi\;\subseteq\;L^2(\mathbb R^n)$ which satisfies
\newline $1.\;\langle {\mathcal R}(h_1), {\mathcal R}(h_2) \rangle_{C(\mathbb T^n)}\;=\;\langle h_1,h_2 \rangle_{C(\mathbb T^n)},$
for all $h_1\;h_2\;\in\;X(q,0,0,\cdots,0,a_j,0,\cdots,0),$
\newline and 
\newline $2.\;{\mathcal R}(X(q,0,0,\cdots,0,a_j,0,\cdots,0)),\;\subseteq\; D {\mathcal R}(X(q,0,0,\cdots,0,a_j,0,\cdots,0)),$ where $D$ is the conjugate of the dilation operator with respect to the dilation matrix $A$ by the Fourier transform. 
\end{theorem} 
\begin{proof}
Without loss of generality we let $j=n-1$ since the proof for the other cases is the same, up to a permutation of the variables.
We recall that the inner product on $\Xi$ is defined by
\[
    \langle \lambda, \mu \rangle_{\mathcal B}(s_1,s_2,\cdots,s_{n-1},t) := \sum_{v \in {\mathbb Z}^n}
                   (\overline{\lambda}\mu)((s_1,s_2,\cdots,s_{n-1},t)-v).
  \]
As in the proof of Theorem 4 of \cite{pacrief2}, we want to find
some appropriate function $\gamma\;\in\;\Xi_n$ and then  define
$${\mathcal R}(h)\;=\;h\;\cdot\;\gamma,\;h\;\in\;X(q,0,0,\cdots,0,a_{n-1}).$$
As in the proof of that theorem , one checks that in order for Equation $1$ in the statement of the theorem to hold, it is necessary to have
 $$\sum_{(m_1,\cdots,m_n) \in {\mathbb Z}^n}|\gamma((s_1-m_1,\cdots,s_{n-1}-m_{n-1},t-m_nq))|^2=1.$$
\newline Then, using the fact that ${\widehat {D_A}}=D_{A^{-1}},$  one checks that if one wants $2$ to hold, 
it is necessary that  
$$D_A({\mathcal R}(X(q,0,0,\cdots,0,a_{n-1}))\;\subseteq\;{\mathcal R}(X(q,0,0,\cdots,0,a_{n-1})),$$
which will as in the proof of Theorem 4 of \cite{pacrief2} will imply the following condition on the function $\gamma:$
$$\{(\prod_{m=1}^n|d_m|)^{1/2}\gamma(d_1s_1,\cdots,d_{n-1}s_{n-1},d_nt)f(d_1s_1,\cdots,d_{n-1}s_{n-1},d_nt):\;f\;\in\;X(q,a)\}\;\subseteq$$
$$\;\{\gamma(s_1,\cdots,s_{n-1},t)f(s_1,\cdots,s_{n-1},t):\;f\;\in\;X(q,0,\cdots,0,a_{n-1})\}.$$
For $f\;\in\;X(q,0,\cdots,0,a_{n-1}),\;f$ must satisfy the equalities
$$f(d_1(s_1-1),\cdots,\;d_{n-1}(s_{n-1}-1),d_nt)\;=\;f(d_1s_1,\cdots,\;d_{n-1}s_{n-1},d_nt),$$ and  
$$f(d_1s,\cdots,d_{n-1}s_{n-1},d_n(t-q))\;=\;e(d_{n-1}d_nqs)f(d_1s,\cdots,d_{n-1}s_{n-1},d_n(t-q)).$$
Once more we check that in order for $$(\prod_{l=1}^n|d_l|)^{1/2}\gamma(d_1s_1,\cdots,d_{n-1}s_{n-1},d_nt)f(d_1s_1,\cdots,d_{n-1}s_{n-1},d_nt)$$  to be an element of ${\mathcal R}(X(q,0,\cdots,0,\;a_{n-1})$ for every $f\;\in\;X(q,0,\cdots,0,\;a_{n-1}),$
it is sufficient that there exist $\mbox{\v m}\;\in\;X(q,0,\cdots,0,\;-(d_{n-1}d_n-1)a_{n-1})$ satisfying 
$$\gamma(d_1s_1,d_2s_2,\cdots,d_{n-1}s_{n-1},d_nt)\;=\;\frac{\mbox{\v m}(s_1,\cdots,s_{n-1},t)}{\sqrt{d}} \gamma(s_1,\cdots,s_{n-1},t),$$
for $d=|d_1d_2\cdots d_n|.$ 
But now we just use the tensor product construction to do this, that is, we find using the methods described in classical wavelet theory and summarized in \cite{pacrief} a continuous function $m'$ of the $n-2$ real variables 
$s_1,\;s_2,\;\cdots,\;s_{n-2},$ which is periodic modulo $\mathbb Z^{n-2}$ and such that $m'$ is continuously differentiable and non-zero in a large enough neighborhood of the origin, $$m'(0,\cdots,0)=\sqrt{|d_1 d_2 \cdots d_{n-2}|},$$
$$\sum_{l=1}^{n-2}\sum_{i_l=0}^{|d_l|-1}|m'(s_1-\frac{i_1}{|d_1|},s_2-\frac{i_2}{|d_2|},\cdots,s_{n-2}-\frac{i_{n-2}}{|d_{n-2}|})|^2\;=\;|d_1 d_2 \cdots d_{n-2}|$$
(we can choose $m'$ to be a tensor product of $n-2$ standard Haar filter functions each with dilation factor $d_i$ constructed in Section 2), 
and define ${\tilde {m}}\;\in\;X(q,-(d_{n-1}d_n-1)a_{n-1}),\; {\tilde {m}}:\;\mathbb R^2\;\rightarrow\;\mathbb C\;$ as in the proof of Theorem 7 of \cite{pacrief2}, with $d_1,d_2,$ and $a$ in that theorem replaced by $d_{n-1},\;d_n,$ and $a_{n-1}$ respectively.   Finally, let $\mbox{\v m}:\;\mathbb R^n\;\rightarrow\;\mathbb C$ be defined by 
$$\mbox{\v m}(s_1,\cdots,s_{n-1},t)\;=\;m'(s_1,\cdots,s_{n-2}){\tilde {m}}(s_{n-1},t).$$
Then if we let
$$\gamma(s_1,\cdots,s_{n-1},t)\;=\;\prod_{i=1}^{\infty}\frac{\mbox{\v m}(s_1/(d_1)^i,\cdots,s_{n-1}/(d_{n-1})^i,t/(d_n)^i)}{\sqrt{d}},$$
it is clear from the tensor product construction that $\gamma$ satisfies the desired conditions, by separation of variable techniques.
\end{proof}
Just as in the $2$-dimensional case, by modifying the proof of Theorem 6 of \cite{pacrief2} we can define a nested sequence of finitely generated projective $C(\mathbb T^n)$-modules
$$\{V_i\}_{i=0}^{\infty}$$
by setting $$V_i\;=\;D^i(V_0).$$
Using the exact same method as in Theorem 6 of \cite{pacrief2} for the two dimensional case, we  associate an ordinary Hilbert space multiresolution analysis to $V_0$ and we obtain the decomposition
$$\Xi\;=\; V_0\oplus_{i\geq 0}\; W_i,$$
where for each $i,\;W_i$ is the finitely generated projective $C(\mathbb T^n)$-module defined by 
$$W_i\;=\;V_{i+1}\;\ominus\; V_i.$$ 
We now want to calculate the image of $V_1\;=\;D(V_0)$ in $K_0(C(\mathbb T^n)),$ where $V_0$ is isomorphic to $X(q,0,0,\cdots,0,a_j,0,\cdots,0)$ as in Theorem \ref{thm quasi3}.  Again without loss of generality, by permuting the variables, we can take $V_0$ isomorphic to $X(q,0,0,\cdots,0,a_{n-1}).$  Somewhat surprising generalizations of Theorem 7 of \cite{pacrief2} now arise:   
\begin{theorem}
\label{thm K0n}
Let $A$ be the $n\times n$ diagonal matrix 
$$A\;=\;\left(\begin{array}{rrrr}
d_1&0&\cdots&0\\
0&d_2&\cdots&0\\
\cdot&\cdot&\cdots&\cdot\\
0&0&\cdots&d_n
\end{array}\right),$$ where $d_1,\;d_2,\cdots,\;d_n\;\in\;\mathbb N$ and
$|d_j|>1,\;1\;\leq\;j\;\leq\;n,$
and let $V_0$ be the image of $X(q,0,0,\cdots,0,a_{n-1})$ in $\Xi\;\subseteq\;L^2(\mathbb R^n)$ defined in Theorem \ref{thm quasi3}, so that $V_1\;=\; D(V_0)\;\supseteq\;V_0.$  Then 
$V_1$ is isomorphic to 
$$X(|\mbox{det}(A)|q,0,0,\cdots,0,\prod_{j=1}^{n-2}|d_j|\cdot \mbox{sign}(d_{n-1}d_n)\;a_{n-1})$$
as a finitely generated projective $C(\mathbb T^n)$ module.
\end{theorem}
\begin{proof}
Let us view $C(\mathbb T^n)$ as a tensor product $C^{\ast}$-algebra,
$$C(\mathbb T^n)\;\cong\;C(\mathbb T^{n-2})\otimes_{\mathbb C}\;C(\mathbb T^2),$$
and we view $X(q,0,0,\cdots,0,a_{n-1})$ as a tensor product of modules with respect to this decomposition:
$$X(q,0,0,\cdots,0,a_{n-1})\;\cong\;C(\mathbb T^{n-2})\otimes_{\mathbb C}\;X(q,a_{n-1}),$$
by recalling that $C(\mathbb T^{n-2})$ is a free module over itself.
With respect to this notation, the module $V_0$ can also be viewed as a tensor product:
$$V_0\;\cong\;V_0'\;\otimes_{\mathbb C}\; V_0''.$$
Here by construction, $V_0'$ is a singly generated free module over $C(\mathbb T^{n-2})$ and can alternatively be viewed as the intersection of $\Xi'\;\subseteq\;L^2(\mathbb R^{n-2})$ with an ordinary multi-resolution analysis level zero space $V_0'$ corresponding to translation by $\mathbb Z^{n-2}$ and dilation by the $(n-2)\;\times\;(n-2)$ diagonal matrix 
$$A'\;=\;\left(\begin{array}{rrrr}
d_1&0&\cdots&0\\
0&d_2&\cdots&0\\
\cdot&\cdot&\cdots&\cdot\\
0&0&\cdots&d_{n-2}
\end{array}\right).$$
The module $V''_0$ corresponds to the image of $X(q,a_{n-1})$ in $\Xi\;\subseteq\;L^2(\mathbb R^2)$ constructed in Theorem 4 of \cite{pacrief2}.
With respect to this decomposition, the unitary operator $D$ can be also be written as a tensor product:
$$D\;=\;D'\;\otimes\;D'',$$
where $D'$ is the dilation operator in the frequency domain corresponding to the matrix $A'$ above, and $D''$ is the dilation operator in the frequency domain corresponding to the $2\times 2$ diagonal matrix
$A''\;=\;\left(\begin{array}{rr}
d_{n-1}&0\\
0&d_n
\end{array}\right).$
Thus we can view $D(V_0)$ as being isomorphic to the tensor product of modules
$$D'(V_0')\;\otimes_{\mathbb C}\;D''(V_0'').$$
Now tensor product techniques and 
standard results from multiresolution analysis show that 
$$D'(V_0')\;\cong\;\oplus_{k=1}^{|\mbox{det}(A')|}[C(\mathbb T^{n-2})]_k$$ as a $C(\mathbb T^{n-2})$ module, 
and the results of Theorem 7 of \cite{pacrief2} show that
$$D''(V_0'')\;\cong\;X(|d_{n-1}d_n|q,\mbox{sign}(d_{n-1}d_n)\;a_{n-1})$$ as a $C(\mathbb T^2)$-module.
Thus $D(V_0)$ is isomorphic to the tensor product of modules
$$\oplus_{k=1}^{|\mbox{det}(A')|}[C(\mathbb T^{n-2})]_k\;\otimes_{\mathbb C}\;X(|d_{n-1}d_n|q,\mbox{sign}(d_{n-1}d_n)\;a_{n-1}).$$
But for each $k,\;1\;\leq\;k\;\leq\;n-2$ 
$$[C(\mathbb T^{n-2})]\;\otimes\;X(|d_{n-1}d_n|q,\mbox{sign}(d_{n-1}d_n)\;a_{n-1})\;\cong\;X(|d_{n-1}d_n|q,0,0,\cdots,0,\mbox{sign}(d_{n-1}d_n)\;a_{n-1})$$ as a $C(\mathbb T^n)$-module.  So by the distributive properties for tensor products and direct sums of $(C(\mathbb T^n)$-modules, 
$$\;(\oplus_{k=1}^{|\mbox{det}(A')|}[C(\mathbb T^{n-2})]_k)\;\otimes_{\mathbb C}\;X(|d_{n-1}d_n|q,\mbox{sign}(d_{n-1}d_n)\;a_{n-1})\;\cong$$
$$\;\oplus_{k=1}^{|\mbox{det}(A')|}[X(|d_{n-1}d_n|q,0,0,\cdots,0,\mbox{sign}(d_{n-1}d_n)\;a_{n-1})]_k\;\cong$$
$$\;X((|\mbox{det}(A)|\cdot q,0,0,\cdots,0,|\mbox{det}(A')|\cdot \mbox{sign}(d_{n-1}d_n)\;a_{n-1}),$$
by an easy generalization of Theorem 6 of \cite{pacrief2} to this situation.
Hence 
$$V_1\;\cong\;X(|\mbox{det}(A)|q,0,0,\cdots,0,|\mbox{det}(A')|\mbox{sign}(d_{n-1}d_n)\;a_{n-1})$$
as finitely generated $C(\mathbb T^n)$-modules, as desired.
\end{proof}
For an example, we consider the case $n=3.$  In this situation, given a non-free $C(\mathbb T^3)$-modules $V_0,$ a non-free module $W_0\;=\; V_1 \;\ominus\; V_0$  can easily arise even when all of the dilation factors have positive signs.
\begin{example}
Let $n=3,$ and $A\;=\;\left(\begin{array}{rrr}
2&0&0\\
0&2&0\\
0&0&2
\end{array}\right).$  For $a\;\in\;\mathbb Z\;\backslash\;\{0\},$ embed $X(1,0,a)$ as a finitely generated projective submodule $V_0$ of $\Xi\;\subseteq\;L^2(\mathbb R^3)$ by means of Theorem \ref{thm quasi3}. 
Then Theorem \ref{thm K0n} shows that $V_1\;=\; D(V_0)$ is isomorphic to $X(8,0,2a),$ so that
$$W_0\;=\; V_1 \;\ominus\; V_0 \;\cong\;X(7,0,a),$$
a finitely generated projective $C(\mathbb T^3)$-module which is not free.
In the last calculation we used the fact that finitely generated projective modules over $C(\mathbb T^n)$ satisfy the cancellation property for $n\;\leq\;4.$
\end{example}
We can now easily identify the isomorphism class of the $V_i$ as $C(\mathbb T^n)$-modules, for $n\geq 2.$ The proof is left to the reader as an exercise:
\begin{corollary}
\label{cor K0n}
Let $q\;\in\;\mathbb N\;,\;a_{n-1}\;\in\;\mathbb Z,$ integers $d_1,d_2,\cdots,d_n,\;A\;=\;
\left(\begin{array}{rrrr}
d_1&0&\cdots&0\\
0&d_2&\cdots&0\\
\cdot&\cdot&\cdots&\cdot\\
0&0&\cdots&d_n
\end{array}\right),$ and $V_0$ 
be as in \ref{thm K0n}, so that
$V_0$ is a submodule of $\;\Xi\;\subseteq\;L^2(\mathbb R^n)$ isomorphic to $X(q,0,0,\cdots,0,a_{n-1}).$ For each $i\;\in\;\mathbb N$ define the finitely generated projective $C(\mathbb T^n)$-module $V_i$ by 
$$V_i\;=\;D^i(V_0).$$
Then $$V_i\;\cong\;X([\prod_{j=1}^n|d_j|]^iq,0,0,\cdots,0,[\prod_{j=1}^{n-2}|d_j|]^i[\mbox{sign}(d_{n-1}d_n)]^i
a_{n-1}),$$
so that 
$$W_i=V_{i+1}\ominus V_i$$
is stably isomorphic to 
$$X([(\prod_{j=1}^n|d_j|)^{i+1}-(\prod_{j=1}^n|d_j|)^{i-1}]q,([\prod_{j=1}^{n-2}|d_j|]^{i+1}[\mbox{sign}(d_{n-1}d_n)]^{i+1}-[\prod_{j=1}^{n-2}|d_j|]^i[\mbox{sign}(d_{n-1}d_n)]^i)a).$$
\end{corollary}
We note that for $n\;\geq\;5,$ in general we can only calculate the {\bf stable} isomorphism type of $W_i\;=\; V_{i+1}\;\ominus\; V_i,$  because finitely generated projective $C(\mathbb T^n)$ modules do not have the cancellation property in general for $n\;\geq\;5.$
We note that this particular exercise would combine the methods of 
Theorem \ref{thm K0n} and  Theorem 8 of \cite{pacrief2} for the calculation.

%%%%%%%%%%%%%%%%%%%%%%%%%%%%%%%%%%%%%%%%%%%%%%%%%%%%%%%%%%%%%%%%%%%%%%%%%%%%%%%
\section{Projective multiresolution analyses for conjugates of diagonal dilation matrices}

Without much extra effort, we can obtain results for dilation matrices that are conjugates of diagonal dilation matrices by elements in $SL(n,\mathbb Z).$  That is, given any $S\;\in\;SL(n,\mathbb Z),$ we can define the $C(\mathbb T^n)$-module
$$X_S(q,0,0,\cdots,a_j,0,\cdots,0)=D_S(X(q,0,0,\cdots,a_j,0,\cdots,0))$$ and inject it as a submodule of $\Xi.$  We are indebted to Dana Williams for the suggestion that this last fact should be true, and for showing us how to prove a major part of it. We first prove the result for $n=2,$ and then move on to the case for general $n.$

\begin{proposition}
\label{prop simdiag}
Let $q\;\in\;\mathbb N$ and $a\;\in\;\mathbb Z$ be fixed, let and $S\;\in\;SL(2,\mathbb Z),$ and set $M\;=\;S^{-1}AS,$ where 
$$A=\;\left(\begin{array}{rr}
d_1&0\\
0&d_2
\end{array}\right),$$
and $d_1,d_2$  are integers both having absolute value strictly greater than $1$. Write $D={\mathcal F} D_A {\mathcal F}^{-1}$ and $D'= {\mathcal F} D_M {\mathcal F}^{-1} .$ Then there exists a finitely generated projective $C(\mathbb T^2)$ submodule of $\Xi$ isomorphic to $X(q,a)$ denoted by  $V_0',$ such that $V_{1}'\;=\;D'(V_0')\;\supseteq\; V_0',$ and more generally, $V_{i}'\;=\;D'(V_{i-1}')\;\supseteq\; V_{i-1}',$ for all $i\;\in\;\mathbb N,$ with 
$$\overline {V_0'\cup\;\cup_{i\in\mathbb N}V_i'}\;=\;\Xi.$$
Furthermore, $$V_{i}'\;\cong\;X(|\text{det}(M)|^iq,[\text{sign(det}(M))]^i\;a),$$ so that 
$$W_0\;=\;V_{1}\ominus V_0\;\cong\;X((|\text{det}(M)|-1)q,(\text{sign}(\text{det}(M))-1)a).$$
\end{proposition}
\begin{proof}
Let the map ${\mathcal R}:\;X(q,a)\;\rightarrow\; V_0$ and the nested sequence of finitely generated projective $C(\mathbb T^2)$-modules 
$\{V_i= D^i(V_0)\}_{i=0}^{\infty}$ be as in 
defined in Theorem 6 of \cite{pacrief2}, so that 
$$\overline{V_0\;\cup\;\cup_{i\in\mathbb N}V_i}\;=\;\Xi.$$
 Note that 
$$D'\;=\;D_{S^t}\circ D \circ D_{(S^t)^{-1}},$$
since ${\mathcal F} D_R {\mathcal F}^{-1}\;=\;D_{(R^t)^{-1}}$ for all $R\;\in\;GL(n,\mathbb R).$
\newline   Define $V_0'\;=\;D_{S^t}(V_0).$
Then 
$$D_{S^t}(V_0)\;=\;\{(\gamma\circ S^t)\cdot (h\circ S^t),\;h\;\in\;X(q,a)\},$$
where $\gamma$ is the function constructed in Theorem 6 of \cite{pacrief2}, so that
$(\gamma\circ S^t)\;\cdot\;(h\circ S^t)\;\in\;\Xi$ for all $h\;\in\;X(q,a).$ Hence  $V_0'\;\subseteq\;\Xi.$
Furthermore, if $(\gamma\circ S^t)\cdot (h\circ S^t)\;\in\; V_0',$
and $g\;\in\;C(\mathbb T^2),$ then
$$[(\gamma\circ S^t)\cdot (h\circ S^t)]g(s,t)\;=\;[(\gamma\circ S^t) (h)\cdot(g\circ (S^t)^{-1})\circ S^t](s,t).$$
Since $g\circ (S^t)^{-1}\;\in\;C(\mathbb T^2),\;(h)\cdot (g\circ (S^t)^{-1})\;\in\;X(q,a),$ so that
$$[(\gamma\circ S^t)\cdot (h\circ S^t)]\cdot g\;\in\;V_0',$$
and $V_0'$ is a $C(\mathbb T^2)$-submodule of $\Xi.$

Now observe that 
$$D'(V_0')\;=\;D_{S^t}\circ D(V_0)\;\supseteq \;D_{S^t}(V_0)\;=\;V_0'.$$

We finally remark that $V_0'\;\cong\;X(q,a).$ 
Indeed for any $R\;\in\;SL(2,\mathbb Z),\;q\;\in\;\mathbb N,$ and $a\;\in\;\mathbb Z,$ consider the following $C(\mathbb T^2)$
module:
$$X_R(q,a)\;=\;\{h\circ\;R:\;h\;\in\;X(q,a)\}.$$
The right action of $C(\mathbb T^2)$ on $X_R(q,a)$ is defined as usual by pointwise multiplication.
One easily checks that $f\cdot g\;\in\;X_R(q,a)$ whenever $f\;\in\;X_R(q,a)$ and $g\;\in\; C(\mathbb T^2),$ 
and another calculation shows that the map on $X_R(q,a)\times X_R(q,a)$ defined by
$$\langle f_1,\;f_2\;\rangle_{C(\mathbb T^2)}(s,t)\;=\;\sum_{j=0}^{q-1}{\overline {f_1(s,t+j)}}f_2(s,t+j)$$
gives a $C(\mathbb T^2)$-valued inner product on $X_R(q,a),$ 
with 
$$\langle f_1,\;f_2\cdot g\;\rangle_{C(\mathbb T^2)}\;=\;\langle f_1,\;f_2\;\rangle_{C(\mathbb T^2)}\cdot\;g,\;\text{for}\;f_1,f_2\;\in\;X_R(q,a)$$
and $g\;\in\;C(\mathbb T^2).$
It is evident from the construction that
in terms of $K_0(C(\mathbb T^2)),$
$$[X_R(q,a)]_{K_0(C(\mathbb T^2))}=(R^{-1})^{\ast}[X(q,a)]_{K_0(C(\mathbb T^2))}=[X(q,a)]_{K_0(C(\mathbb T^2))},$$
by Equation \ref{eq act}.
We now define a map ${\mathcal {R'}}:\;X_{S^t}(q,a)\;\rightarrow\;V_0'$ by 
$${\mathcal {R'}}(f)(s,t)\;=\;\gamma\circ S^t(s,t)\cdot f(s,t),$$
where again $\gamma$ is the function constructed in Theorem 6 of \cite{pacrief2}. It is clear that ${\mathcal {R'}}$ preserves the module structure, and we can check that ${\mathcal {R'}}$
preserves the $C(\mathbb T^2)$-valued inner product as well.  Let $f_i\;=\;h_i\circ\;S^t,$ where $h_i\;\in\;X(q,a),\;i\;=\;1,2.$ Then 
$$\langle f_1,\;f_2\;\rangle_{C(\mathbb T^2)}\;=\;\langle h_1,\;h_2\;\rangle_{C(\mathbb T^2)}\circ\;S^t,$$ 
and a similar calculation shows that
$${\mathcal {R'}}(f_i)\;=\;{\mathcal {R}}(h_i)\circ\;S^t,\;i\;=\;1,2.$$
Hence, for $f_1,\;f_2\;\in\;X_{S^t}(q,a)$ and $h_1,\;h_2\;\in\;X(q,a)$ as above, 
$$\langle f_1,\;f_2\;\rangle_{C(\mathbb T^2)}(s,t)\;=\;\langle h_1,\;h_2\;\rangle_{C(\mathbb T^2)}\circ\;S^t\;(s,t)$$
$$\;=\;\langle {\mathcal {R}}(h_1), {\mathcal {R}}(h_2) \rangle_{C(\mathbb T^2)}\circ\;S^t\;(s,t)$$
$$\;=\;\sum_{(m,n) \in \mathbb Z^2}{\overline{[(\gamma\circ S^t)\cdot (f_1)]}}\cdot[(\gamma\circ\;S^T)\cdot (f_2)]((s,t)-(m,n))$$
$$\;=\;\langle {\mathcal {R'}}(f_1),\;{\mathcal {R'}}(f_2)\;\rangle_{C(\mathbb T^2)},$$
so that ${\mathcal {R'}}$ preserves the $C(\mathbb T^2)$-valued inner products, as desired.
Since $V_0'$ is isomorphic to $X_{S^t}(q,a)$ and $X_{S^t}(q,a)$ is isomorphic to $X(q,a),$ it follows that $V_0'$ is isomorphic to $X(q,a)$ as a $C(\mathbb T^2)$-module.
\newline Finally, we will show that
$$V_{i}'\;\cong\;X(|\mbox{det}(M)|^iq,[\mbox{sign(det)}(M)]^i\;a).$$
Note that $\mbox{det}(M)\;=\;\mbox{det}(A),$ as $M$ and $A$ are similar, so it is enough to show that
$$V_{i}'\;\cong\;X(|\mbox{det}(A)|^nq,[\mbox{sign(det)}(A)]^i\;a).$$
We do this for the case $i=1,$ and note that higher dimensional cases are similar.  It is enough to show that  $$D_{S^t}\circ D'(V_0)\;\cong\;X(|\mbox{det}(A)|^q,[\mbox{sign(det)}(A)]\;a).$$
For ease of notation, we now write $d_1=b_1|d_1|,\;d_2=b_2|d_2|,$ where $b_1,\;b_2,\;\in\;\{-1,1\}.$  Then $\text{sign[det}(A)]\;=\;b_1b_2.$ 
Let $(d_1,a)=w,$ with $w>0,$ so that we can write $a=r\cdot w,$ and $d_1=y\cdot w,$ where $(r,y)=1,$ i.e. $r$ and $y$ are relatively prime.
Then $|d_1|\;=\;|y\cdot w|\;=\;|y|\cdot w\;=\;yb_1\cdot w,$ so that $|y|\;=\;y\cdot b_1$
Arguments similar to those used in the proof of Theorem 7 of \cite{pacrief2} show that 
$D(V_0)\;\cong\;\oplus_{j=0}^{w-1}[X(|yd_2|q,b_1b_2r)]_j,$
so that 
$$D_{S^t}\circ D(V_0)\;\cong\;\oplus_{k=0}^{w-1}[X_{S^T}(|yd_2|q,b_1b_2r)]_j].$$
Again, the discussion given in the first part of the proof shows that
$$X_{S^t}(|yd_2|q,b_1b_2r)]_j\;\cong\;X(|yd_2|q,b_1b_2r)]_j,$$
so that  
$$D_{S^t}\circ D (V_0)\;\cong\;\oplus_{j=0}^{w-1}[X(|yd_2|q,b_1b_2r)]_j],$$
and since $$\oplus_{j=0}^{w-1}[X(|yd_2|q,b_1b_2r)]_j\;\cong\;X(|d_1d_2|q,\mbox{sign(det}(A))a),$$ the desired result has been proved.
\end{proof}

\begin{theorem}
\label{thm simdiag}
Let $q\;\in\;\mathbb N$ and $a\;\in\;\mathbb Z$ be fixed, and let $M\;\in\;M(n,\mathbb Z)$ and $S\;\in\;SL(2,\mathbb Z)$ be such that
$$SMS^{-1}\;=\;A\;=\;\left(\begin{array}{rrrr}
d_1&0&\cdots&0\\
0&d_2&\cdots&0\\
\cdot&\cdot&\cdots&\cdot\\
0&0&\cdots&d_n
\end{array}\right),$$ where $d_1,\;d_2,\cdots,\;d_n\;\in\;\mathbb N$ and
$|d_j|>1,\;1\;\leq\;j\;\leq\;n.$
Then there exists a finitely generated projective $C(\mathbb T^n)$ submodule of $\Xi$ isomorphic to $X(q,0,0,\cdots,a_j,0,\cdots,0),$ denoted by  $V_0',$ such that $V_1'\;=\;D'(V_0')\;\supseteq\;V_0'$ and  
$$\overline {V_0'\cup\;\cup_{i\in\mathbb N} V_i'}\;=\;\Xi.$$
\end{theorem}

\begin{proof}
Let the $C(\mathbb T^n)$-module $V_0$ which is isomorphic to $X(q,0,0,\cdots,a_j,0,\cdots,0)$ be as in the proof of Theorem \ref{thm quasi3}, and let $D'={\mathcal F}\circ D_M \circ {\mathcal F}^{-1}.$ 
As in the proof of Proposition \ref{prop simdiag} we define 
the $C(\mathbb T^n)$-submodule $V_0'$ which is isomorphic to 
$X_{S^t}(q,0,0,\cdots,a_j,0,\cdots,0)$ and which satisfies the condition
$$V_0'\;\subseteq\;D'(V_0'),$$ where as before $V_0'\;=\;D_{S^t}(V_0).$
The proof is then very similar to that given in the proof of Proposition \ref{prop simdiag}, and we leave details to the reader.
\end{proof}

\begin{remark}
\label{remark higherdim}
Let $S\in SL(2,\mathbb Z),\;n\geq 3,$ and put $B=S^t.$ The main difference in passing to the case $n\geq 3$ is that $X_{B}(q,0,0,\cdots,a_j,0,\cdots,0)$ will not be isomorphic to $X(q,0,0,\cdots,a_j,0,\cdots,0)$ in general, unlike the two-dimensional case.  Recall that 
\newline $[X(q,0,0,\cdots,0,a_j,0,\cdots,0)]_{K_0(C(\mathbb T^n))}$ can be identified with the element $q\;-\;a_je_j\wedge\;e_n\;\in\;\bigwedge^{\mbox{\scriptsize even}}_{\mbox{\scriptsize \bf Z}}
\{e_1,...,e_n\}.$  On the other hand, suppose that we write 
$$B\;=\;(b_{jk})_{1\leq\;j,k\;\leq n}.$$ 
Upon modifying the methods of Proposition 3.4 of \cite{pacrief2}, we see that 
\newline $[X_S(q,0,0,\cdots,a_j,0,\cdots,0)]_{K_0(C(\mathbb T^n))}$
is identified with the element 
\newline $q\;-\;a_j(\sum_{l=1}^n b_{jl}e_l)\wedge(\sum_{k=1}^n b_{nk}e_k)\;\in\;\bigwedge^{\mbox{\scriptsize even}}_{\mbox{\scriptsize \bf Z}}
\{e_1,...,e_n\},$ which of course need not be equal to $q\;-\;a_je_j\wedge\;e_n$ for $n\;\geq\;3.$  
\end{remark}

Suppose further that
$d_1=d_2=\cdots=d_{n-1}=d_n.$  Then $A$ commutes with any matrix $B\;\in\;SL(n,\mathbb Z),$ i.e. $B^{-1}AB=A.$ In this way we can obtain a large variety of $C(\mathbb T^n)$ submodules $V_0'$ of $\Xi$ such that 
$$V_0'\;\subseteq\;D(V_0').$$ 

We develop this idea for the case $n=3.$
\begin{example}
\label{ex all3}
Let $n=3$ and let our dilation matrix be the standard dilation by $d$ in all variables, i.e. $A\;=\;\left(\begin{array}{rrr}
d&0&0\\
0&d&0\\
0&0&d
\end{array}\right).$  Let $B\;=\;(b_{jk})_{1\leq\;j,k\;\leq 3}$ be an arbitrary element of $SL(3,\mathbb Z).$    
Then the above argument shows that we can take 
$$V_0'\;\cong\;X_B(q,0,1)$$ for any $q\;\in\;\mathbb N,$ and build a projective multiresolution analysis for dilation by $A.$
One calculates from the formulas given above that
$$[X_B(q,0,a)]_{K_0(C(\mathbb T^3))}\;=\;[q\;+\;a[b_{22}b_{31}-b_{21}b_{32}]e_1\wedge e_2\;+\;a[b_{23}b_{31}-b_{21}b_{33}]e_1\wedge e_3$$
$$\;+\;a[b_{23}b_{32}-b_{22}b_{33}]e_2\wedge e_3]_{K_0(C(\mathbb T^3))}.$$
\end{example}
With this example in hand, whenever the $3\times 3$ dilation matrix is a constant integral multiple of the identity, we can easily embed projective modules isomorphic to any nonzero representative of the positive cone of $[K_0(C(\mathbb T^3)]]$ as the initial space of one of our projective  multiresolution analyses, as long as we keep in mind that every nonzero member of this positive cone can be represented as $[q+c_1e_1\wedge e_2+c_2e_1\wedge e_3+c_3e_2\wedge e_3]$ for $q\;\in\;\mathbb N$ and $c_1,c_2,c_3\;\in\;\mathbb Z$ via its exterior algebra representation in $\bigwedge^{\mbox{\scriptsize even}}_{\mbox{\scriptsize \bf Z}}
\{e_1,e_2,e_3\}.$

\begin{theorem}
\label{thm all3}
Let $n=3$ and let our dilation matrix be the standard dilation by $d$ in all variables, i.e. $A\;=\;\left(\begin{array}{rrr}
d&0&0\\
0&d&0\\
0&0&d
\end{array}\right),$ where d is an integer whose absolute value is greater than $1.$  Then, identifying $K_0(C(\mathbb T^3))$ with $\bigwedge^{\mbox{\scriptsize even}}_{\mbox{\scriptsize \bf Z}}
\{e_1,e_2,e_3\}\;=\;\bigwedge^0_{\mbox{\scriptsize \bf Z}}\{e_1,e_2,e_3\}\oplus \bigwedge^2_{\mbox{\scriptsize \bf Z}}\{e_1,e_2,e_3\},$ given any $[q+c_1e_1\wedge e_2+c_2e_1\wedge e_3+c_3e_2\wedge e_3]\;\in\;[K_0(C(\mathbb T^3))]^+,$ there exists a finitely generated projective $C(\mathbb T^3)$-module ${\widetilde {V_0'}}\;\subseteq\;\Xi$ with  
$$V_0'\;\cong\;[q+c_1e_1\wedge e_2+c_2e_1\wedge e_3+c_3e_2\wedge e_3]$$
such that
$$D(V_0')\;\supseteq\; V_0'.$$ Hence, setting 
$$V_n'\;=\; D^n(V_0')$$
for $n\;\geq\;0,$ the nested sequence $\{V_n'\}$ gives a projective multi-resolution analysis corresponding to the dilation matrix $A.$
\end{theorem}
\begin{proof}
Let the greatest common divisor of $c_1,\;c_2$ and $c_3$ be $a\;\in\;\mathbb N,$ and write $c_1=ax,\;c_2=ay,$ and $c_3=az,$ where $x,\;y,$ and $z$ are relatively prime.  By Example \ref{ex all3}, in order to prove the theorem, it is sufficient to find $B\;=\;(b_{jk})_{1\leq\;j,k\;\leq 3}\;\in\;SL(3,\mathbb Z),$ such that
$[b_{22}b_{31}-b_{21}b_{32}]=x,\;[b_{23}b_{31}-b_{21}b_{33}]=y,$ and 
$[b_{23}b_{32}-b_{22}b_{33}]=z.$  The three terms involving the $b_{ij}$ are nothing more than $-1$ times the cofactors of $B$ obtained from its last two rows, so the question then becomes one of whether it is possible to find 
$B\;\in\;SL(3,\mathbb Z)$ with specified relatively prime cofactors along the last two rows.  Of course this is possible; first choose $b_{11},b_{12},$ and $b_{13}\;\in\;\mathbb Z$ such that 
$$-b_{11}z+b_{12}y-b_{13}x\;=1.$$
Now let $(x,z)\;=\nu,$ so that $(\nu,y)=1$ as $x,y,$ and $z$ are relatively prime.  Write $x\;=\;M\alpha$ and $z\;=\;\nu\beta$ where $(\alpha,\beta)=1.$  Since $\alpha$ and $\beta$ are relatively prime, we can find $\sigma,\;\tau\;\in\;\mathbb Z$ such that
$\alpha\tau+\beta\tau\;=\;y$.
Now we set
$$b_{21}=-\alpha,\;b_{22}\;=\;0,\;b_{23}=\beta,$$
$$b_{31}=\sigma,\;b_{32}=\nu,\;b_{33}=\tau.$$
Then one checks that $[b_{22}b_{31}-b_{21}b_{32}]=x,\;[b_{23}b_{31}-b_{21}b_{33}]=y,\; 
[b_{23}b_{32}-b_{22}b_{33}]=z,$ and that the matrix $B$ with entries as defined above satisfies 
$$\mbox{det}(B)\;=\;1,$$ so is an element of $SL(3,\mathbb Z),$ as desired.
Finally, setting $V_0'\;=\;X_B(q,0,a),$
we have 
$$D(V_0')\;\supseteq\; V_0'$$
and
$$[V_0']_{K_0(C(\mathbb T^3))}\;=\;[q\;+\;a[b_{22}b_{31}-b_{21}b_{32}]e_1\wedge e_2\;+\;a[b_{23}b_{31}-b_{21}b_{33}]e_1\wedge e_3$$
$$\;+\;a[b_{23}b_{32}-b_{22}b_{33}]e_2\wedge e_3]_{K_0(C(\mathbb T^3))}$$
$$\;=\;[q+c_1e_1\wedge e_2+c_2e_1\wedge e_3+c_3e_2\wedge e_3]_{K_0(C(\mathbb T^3))},$$
as desired.
We leave it to the reader as an exercise to calculate the isomorphism class
of each $V_n'$ in $K_0(C(\mathbb T^3))$ by using the combined methods of Proposition \ref{thm simdiag} and Corollary \ref{cor K0n}.
\end{proof} 
By rearranging the variables, the above discussion generalizes slightly to allow us to deal with diagonal dilation matrices and finitely generated projective modules whose Chern characters take on the values $q_j-ae_{i_j}\wedge\;e_{k_j}$ for $i_j\;\not=\;k_j$ in the exterior algebra formulation of $K_0(C(\mathbb T^n)).$  Using separation of variable techniques, we can generalize even further as follows: suppose that $[P]\;\in\;K_0(C(\mathbb T^n))^+$ has its representation in $\bigwedge^{\mbox{\scriptsize even}}_{\mbox{\scriptsize \bf Z}}
\{e_1,...,e_n\}$ given by
$$[P]\;=\bigwedge_{j=1}^{l}\;(q_j-a_1e_{i_j}\wedge\;e_{k_j}),$$
where $q_j\;\in\;\mathbb N,\;2l\;\leq\;n,$ and $\{i_1,k_1,i_2,k_2,\cdots,i_l,k_l\}$ are distinct elements $\{1,2,\cdots,n\}.$  Without loss of generality by rearranging the variables we assume that $i_1<k_1<i_2<k_2<\cdots<i_l<k_l.$
Then by repeated use of Theorem \ref{thm K0n} and tensor products, it is possible to find a finitely generated 
projective $C(\mathbb T^n)$-submodule $V_0\;\subseteq\;\Xi$ such that
$$[V_0]\;=\;[P]\;\in\;K_0(C(\mathbb T^n))$$
and $$V_0\;\subseteq\; D(V_0),$$
where $D$ is dilation in the frequency domain by an $n\times n$ diagonal matrix $A$ with integer entries all having absolute value greater than $1.$ We leave details to the reader.

%%%%%%%%%%%%%%%%%%%%%%%%%%%%%%%%%%%%%%%%%%%%%%%%%%%%%%%%%%%%%%%%%%%%%%%%%%%%%%%
\section{Building module frames from projective multi-resolution analyses}

In this section, we present a useful algorithm for constructing a module frame for $\Xi$ from a projective multiresolution analysis corresponding to a dilation matrix $A.$ 

Before beginning, we fix a bit of notation:  if $v\in \mathbb Z^n,$ recall that the translation operator $T_v$ on $L^2(\mathbb R^n)$ is defined by 
$$T_v(f)(x)\;=f(x-v),\;f\in L^2(\mathbb R^n).$$ 
For any $v\in \mathbb Z^n,$ let $\varepsilon_v$ denote the Fourier transformed version of $T_v$ acting on $L^2(\mathbb R^n)$ defined by $\varepsilon_v\;=\;{\mathcal F}\circ T_v\circ {\mathcal F}^{\ast}.$  A standard calculation shows that $\varepsilon_v$ is the multiplication operator on $L^2(\mathbb R^n)$ defined by 
$$\varepsilon_v (\xi) (t)\;=\;e(-v\cdot t)\xi(t),\;\xi\in L^2(\mathbb R^n),$$ where recall 
$e$ is the exponential function defined on $\mathbb R$ by $e(r)=e^{2\pi ir}.$

We now present the theorem on the construction of module frames. 
\begin{theorem}
\label{thm pmramodfram}
Let $A\;\in\;GL(n,\mathbb Q)$ be a $n\times n$ dilation matrix with integer entries, let  $|\text{det}(A)|=d,$ and suppose that  $\{V_j\}_{j=0}^{\infty}$ is a projective multiresolution analysis for $\Xi$ corresponding to dilation by $A.$  Suppose that  $\{\Phi_1,\Phi_2,\cdots,\Phi_{s}\}$ is a module frame for $V_0$, and that $\{\Psi_{1},\Psi_{2},\cdots,\Psi_{r}\}$ is a module frame for the finitely generated projective $C(\mathbb T^n)$-module $W_0\;=\; V_1\ominus\; V_0.$ Then it is possible to enumerate a countably infinite subset $\{v_l:\;l\geq 0\}\;\subseteq\;\mathbb Z^n$ such that for each $i\in\;\mathbb N\cup\{0\},$ $\{v_l:\;0\;\leq\;l\;\leq\;d^i-1\}$ is a set of coset representatives for $\mathbb Z^n/A^i(\mathbb Z^n).$ 
Furthermore, letting $D$ be dilation by $A$ in the frequency domain, for each $i\in\;\mathbb N,$
$$\{D^{i}\varepsilon_{v_{k_1}}(\Psi_1),D^{i}\varepsilon_{v_{k_2}}(\Psi_2),\cdots, D^{i}\varepsilon_{v_{k_r}}(\Psi_r):\;0\leq k_l\leq d^i-1,\;1\leq l\leq r\}$$
forms a module frame for the finitely generated projective $C(\mathbb T^n)$-module $W_i\;=\; V_{i+1}\ominus\; V_i,$ so that a module frame for all of 
$\Xi$ is given by the set
$$\{\Phi_1,\Phi_2,\cdots,\Phi_{s}\}\;\cup\;\cup_{j=0}^{\infty}\cup_{k=1}^r\{D^j\varepsilon_{v_l}(\Psi_k):\;0 \leq\;l \;\leq \;d^j-1 \}.$$
\end{theorem}
\begin{proof}
We start the proof by showing how to construct the appropriate sets of coset representatives, enumerated in a consistent fashion.  Let $\{\beta_i:\;0\leq i\leq d-1\}$ be a set of coset representatives for $\mathbb Z^n/A(\mathbb Z^n)$ with $\beta_0\;=\;0.$
Fix $b\;\in\;\mathbb N.$ Then it is easy to check that $\{\sum_{i=0}^{b-1}A^i(\beta_{b_i}):\;0\leq b_i\leq d-1\}$ gives a family of coset representatives for
$\mathbb Z^n/A^b(\mathbb Z^n).$ The set 
$\{v_l:\;0 \;\leq\;l\;\leq\;d^i-1\}$ is given by giving $l$ its $d$-adic expansion in terms of powers of the matrix $A$ and the coset representatives $\beta_i,$ i.e. if we write $l=\sum_{i=0}^{b-1}(b_i)d^i,$ where $b_i\;\in\;\{0,1,\cdots,d-1\},$ we then set
$v_l\;=\;\sum_{i=0}^{b-1}A^i(\beta_{b_i}).$  One can check that as one passes from the enumeration of coset representatives of $\mathbb Z^n/A^b(\mathbb Z^n)$ to the enumeration of representatives for $\mathbb Z^n/A^{b+1}(\mathbb Z^n),$ this definition of $v_l$ is consistent, i.e. it does not depend on $b.$
 
We will establish the first claim of the theorem by induction.  The statement is true for $i=0,$ by hypothesis. We prove it is true for $i=1,$ as the proof of this fact will serve as a model for the general case.  We note that $V_{2}\;=\; V_{1}\oplus\; W_{1}.$  But also $V_{2}\;=\;D(V_{1})\;=\;D(V_{0}\oplus\; W_{0}).$  One easily checks that
$$D(V_{0}\oplus\; W_{0})\;=\; D(V_{0})\oplus\; D(W_{0})\;=\; V_{1}\oplus\; D(W_{0}).$$
We note that  $$W_{1}\;=\;D(W_{0}),$$
since the orthogonal complement of $V_{1}$ in $V_{2}$ with respect to the Hilbert module inner product must be unique.  We know that $\{\Psi_k:\;1\leq k\leq r\}$ is a  module frame for $W_{0}$ and we want to show that
$$\{ D \varepsilon_{\beta_{k_1}}(\Psi_1), D \varepsilon_{v_{k_2}}(\Psi_2),\cdots,  D \varepsilon_{\beta_{k_r}}(\Psi_r):\;0\leq k_l\leq d-1,\;1\leq l\leq r\}$$
is a module frame for $V_{1}.$  Fix $\zeta\;\in\;W_{1}\;=\; D(W_{0}).$
We must show that $$\zeta(x)\;=\;\sum_{k=1}^r\sum_{l=d^j}^{d^{j+1}-1} D \varepsilon_{v_l}(\Psi_k)(x)\langle D \varepsilon_{\beta_l}(\Psi_k), \zeta \rangle_{C(\mathbb T^n)}(x),$$
for $\{\beta_l:\;0\;\leq\;l\;\leq\;d-1\}$ a set of coset representatives for $\mathbb Z^n/A(\mathbb Z^n),$ where without loss of generality we choose $\beta_0=0.$ Since $\zeta\;\in\;D(W_{0}),$ we can write $\zeta\;=\;D(\xi)$ for some $\xi\in\; W_{0}.$
Using the definition of module frame for $W_0,$ we get 
$$\xi(x)\;=\;\sum_{k=1}^r\Psi_k(x)\langle \Psi_k, \xi \rangle_{C(\mathbb T^n)}(x).$$
Since $\zeta(x)=\frac{1}{\sqrt{d}}\xi([A^t]^{-1}(x)),$
to obtain the desired result it is enough to show that
$$\sum_{k=1}^r\sum_{l=0}^{d-1} D \varepsilon_{v_l}(\Psi_k)(x)\langle D \varepsilon_{{\beta_l}}(\Psi_k), \zeta \rangle_{C(\mathbb T^n)}(x)=\;\frac{1}{\sqrt{d}}\sum_{k=1}^r\Psi_k([A^t]^{-1}(x))\langle \Psi_k, \xi \rangle_{C(\mathbb T^n)}([A^t]^{-1}(x)).$$

We now calculate that 
$$\sum_{k=1}^r\sum_{l=0}^{d-1} D \varepsilon_{\beta_l}(\Psi_k)\langle D \varepsilon_{\beta_l}(\Psi_k), \zeta \rangle_{C(\mathbb T^n)}(x)$$
$$=\;\sum_{k=1}^r\sum_{l=0}^{d-1}D \varepsilon_{\beta_l}(\Psi_k)\langle D\varepsilon_{\beta_l}(\Psi_k), \zeta \rangle_{C(\mathbb T^n)}(x)$$
$$=\;\sum_{k=1}^r\sum_{l=0}^{d-1}\varepsilon_{A^{-1}(\beta_l)}D(\Psi_k)(x)\langle \varepsilon_{A^{-1}(\beta_l)}D^{1}(\Psi_k), \frac{1}{\sqrt{d}}\xi\circ [A^t]^{-1} \rangle_{C(\mathbb T^n)}(x)$$
$$=\frac{1}{\sqrt{d^{2}}}\sum_{k=1}^r\sum_{l=0}^{d-1}e(-A^{-1}(\beta_l)\cdot x)(\Psi_k)([A^t]^{-1}(x))\langle \varepsilon_{A^{-1}(\beta_l)}D^{j+1}(\Psi_k),\xi\circ [A^t]^{-1} \rangle_{C(\mathbb T^n)}(x)$$
$$=\frac{1}{d\sqrt{d}}\sum_{k=1}^r\sum_{l=0}^{d-1}e(- A^{-1}(\beta_l)\cdot x)\Psi_k\circ [A^t]^{-1}(x)\cdot$$
$$\cdot\sum_{v\in\;\mathbb Z^n}{\overline{e(- A^{-1}(
\beta_l)\cdot(x+v))\Psi_k\circ [A^t]^{-1}(x+v)}} \cdot\xi([A^t]^{-1}(x+v))$$
$$=\frac{1}{d\sqrt{d}}\sum_{k=1}^r\sum_{l=0}^{d-1}\sum_{v\in\;\mathbb Z^n}\sum_{m=0}^{d-1}e(A^{-1}(\beta_l)\cdot [A^t](v)+\gamma_m)\Psi_k\circ [A^t]^{-1}(x)$$
$$\cdot{\overline{(\Psi_k\circ [A^t]^{-1}(x+[A^t](v)+\gamma_m)}} \cdot\xi([A^t]^{-1}(x+[A^t](v)+v_m)),$$
(where $\{\gamma_m:\;0\leq\;m\leq\;d-1\}$ is a set of coset representatives for $\mathbb Z^n/[A^t](\mathbb Z^n)$ with $\gamma_0=0$)
$$=\frac{1}{d\sqrt{d}}\sum_{k=1}^r\sum_{l=0}^{d^{j+1}-1}\sum_{m=0}^{d^{j+1}-1}e(A^{-1}(\beta_l)\cdot \gamma_m)\Psi_k( [A^t]^{-1}(x)$$
$$\cdot\sum_{v\in\;\mathbb Z^n}\sum_{i=0}^{d-1}{\overline {\Psi_k(([A^t]^{1}(x+\gamma_m)+v))}}\cdot\xi([A^t]^{-1}(x)+v+\gamma_m).$$
By the law of characters on the finite abelian group $\mathbb Z^n/[A^t(\mathbb Z^n)],$ $$\sum_{l=0}^{d^{(j+1)}-1}e((A^{-1}\beta_l)\cdot \gamma_m)=0$$ 
unless $\gamma_m=0,$ that is, unless $m=0,$ and equals $d$ when $m=0,$ in which case $\beta_0=\gamma_0=\vec{0}\in\mathbb Z^n.$
Our sum therefore becomes
$$=\frac{1}{d\sqrt{d}}\cdot d\sum_{k=1}^r\sum_{v\in\;\mathbb Z^n}\Psi_k( [A^t]^{-1}(x){\overline {\Psi_k(([A^t]^{-1}(x)+v))}}\cdot\xi([A^t]^{-1}(x)+v)$$
$$=\;\frac{1}{\sqrt{d}}\sum_{k=1}^r \Psi_k([A^t]^{-1}(x))\langle \Psi_k, \xi \rangle_{C(\mathbb T^n)}([A^t]^{-1}(x))$$
$$=\;\frac{1}{\sqrt{d}}\xi([A^t]^{-1}(x))\;=\; D(\xi)(x)\;=\;\zeta(x),$$
as desired. 

We now suppose that $\{\Psi_{j,1},\Psi_{j,2},\cdots,\Psi_{j,r_j}\}$ is a module frame for 
$W_j,$ the orthogonal complement of $V_j$ in $V_{j+1}$ with respect to the Hilbert module $C(\mathbb T^n_)$-valued inner product. 
The same argument used to show that $W_1= D(W_0)$
shows that $W_{j+1}= D(W_j).$
It is easy now to see that 
$$\{D \varepsilon_{\beta_{k_1}}(\Psi_{j,1}), D \varepsilon_{v_{k_2}}(\Psi_{j,2}),\cdots, D \varepsilon_{\beta_{k_r}}(\Psi_{j,r_j}):\;0\leq k_l\leq d-1,\;1\leq l\leq r_j\}$$
is a module frame for $V_{j+1}.$  This is because the same argument as for the case $j=0$ shows that if  $\zeta\;\in\; W_{j+1}\;=\;D(W_{j}),$ then
$\zeta=D(\xi)$ for some $\xi\in W_j,$
and 
$$\sum_{k=1}^r\sum_{l=0}^{d-1} D \varepsilon_{\beta_l}(\Psi_{j,k})(x)\langle D \varepsilon_{{\beta_l}}(\Psi_{j,k}), \zeta \rangle_{C(\mathbb T^n)}(x)=\;\frac{1}{\sqrt{d}}\sum_{k=1}^{r_j}\Psi_{j,k}([A^t]^{-1}(x))\langle \Psi_{j,k}, \xi \rangle_{C(\mathbb T^n)}([A^t]^{-1}(x)),$$ from which follows the desired statement about the module frame for $W_{j+1}.$  

We now prove by induction that $\{D^{i} \varepsilon_{v_l}(\Psi_k)\;0\leq l\leq d^i-1,\;1\leq k\leq r\}$ forms a module frame for the finitely generated projective $C(\mathbb T^n)$-module $W_i,$ where here $\{v_l:\;0\;\leq\;l\;\leq\;d^i-1\}$ is a set of coset representatives for $\mathbb Z^n/A^i(\mathbb Z^n)$ as defined in the first paragraph of the proof of the theorem. The statement is true for $i=0,$ by hypothesis, and we now assume it is true for $i=j,$ so that $\{D^{j} \varepsilon_{v_l}(\Psi_k)\;0\leq l\leq d^j-1,\;1\leq k\leq r\}$ forms a module frame for the finitely generated projective $C(\mathbb T^n)$-module $W_j,$ where here $\{v_l:\;0\;\leq\;l\;\leq\;d^j-1\}$ is a set of coset representatives for $\mathbb Z^n/A^j(\mathbb Z^n).$  From our remarks above, we know that 
a module frame for $W_{j+1}$ is given by the set 
$$\{D \varepsilon_{\beta_{m}}(D^{j} \varepsilon_{v_l}(\Psi_k)):\;0\leq m\leq d-1,\;0\leq l\leq d^j-1,\;1 \leq k\leq r\}.$$
We now show that this takes on the desired form.
We first note the relations
$$\varepsilon_{\beta_{m}}D^j\;=\; D^j \varepsilon_{A^j(\beta_{m})},,\;\text{and}\; \varepsilon_{A^j(\beta_{m})} \varepsilon_{v_l}\;=\; \varepsilon_{[A^j(\beta_{m})+v_l]}.$$
Thus we see that 
$$\{D^{j+1} \varepsilon_{[A^j(\beta_{m})+v_l]}(\Psi_k)):\;0\leq m\leq d-1,\;0\leq l\leq d^j-1,\;1 \leq k\leq r\}$$
forms a module frame for $W_{j+1}.$
We finally note that as $m$ runs through all the values between $0$ and $d-1,$ and $l$ goes from 
$0$ to $d^j-1,$ the set 
$\{A^j(\beta_{m})+v_l\}$ runs through all the coset representatives of $\mathbb Z^n/A^{j+1}(\mathbb Z^n),$ as explained in our consistent enumeration of coset representatives given in the first paragraph of the proof.
Thus 
$$\{A^j(\beta_{m})+v_l:\;0\leq m\leq d-1,\;0\leq l\leq d^j-1\}\;=\;\{v_l:\;0\leq l\leq d^{j+1}-1\},$$
so that 
$$\{D^{j+1} \varepsilon_{v_l}(\Psi_k)):\;0\leq l\leq d^{j+1}-1,\;1 \leq k\leq r\}$$
forms a module frame for $W_{j+1}.$
The induction step of the first statement of the theorem is done, so we have proved the first statement constructing module frames for $W_i$ for each $i\in\mathbb N\cup\{0\}.$
The rest of the theorem follows from the fact that we have a projective multiresolution analysis for $\Xi,$ so that once we know module frames for $V_0$ and each $W_i,$ the union of these module frames gives a standard module frame for all of $\Xi.$

Finally, we remark that if $W_0$ is a free $C(\mathbb T^n)$ module of rank $r,$ then $W_i$ will be free of rank $rd^i.$  This follows by choosing the module frame for $W_0$ in such a way that $\langle \Psi_i,\Psi_j \rangle_{C(\mathbb T^2)}\;=\;\delta_{i,j}1_{C(\mathbb T^2)},$ performing a calculation involving the law of characters to show that 
$$\langle D \varepsilon_{v_l}(\Psi_i), D \varepsilon_{v_m}(\Psi_j)\rangle_{C(\mathbb T^2)}\;=\;\delta_{i,j}\delta_{l,m}1_{C(\mathbb T^n)},$$
so that $W_1$ is free of rank $rd.$ Using induction exactly as before, we see that $W_i$ is free of rank $rd^i$ for all $i\in \mathbb N.$ 
\end{proof}

%%%%%%%%%%%%%%%%%%%%%%%%%%%%%%%%%%%%%%%%%%%%%%%%%%%%%%%%%%%%%%%%%%%%%%%%%%%%%%%
\section{Conclusion}

The methods we have used so far may seem fairly specialized.
The question then arises as to how to construct examples beyond the special cases of diagonal matrices with integer entries, and their conjugates by elements of $SL(n,\mathbb Z),$ and various special ad-hoc dilation matrices, such as $Q=\left(\begin{array}{rr}
1&1\\
-1&1
\end{array}\right),$ the matrix associated to the quincunx lattice. In the quincunx lattice case, we have been able to construct an initial module $V_0$ in $\Xi$ isomorphic to the projective $C(\mathbb T^2)$ module $X(1,1)$ and this leads to a wavelet module $W_0$ that is a free $C(\mathbb T^2)$-module of rank one. Letting $\phi_1$ and $\phi_2$ be the module frame for $V_0$ as constructed following the proof of Proposition 17 in \cite{pacrief2}, and letting $\{\psi\}$ be a module frame for the singly-generated module $W_0,$  Theorem \ref{thm pmramodfram} gives a way to construct a module frame of $\Xi$ that corresponds to the projective multiresolution analysis.

  By using recent constructions of smooth scaling functions in the ordinary multiresolution analysis setting by  M. Bownik \cite{Bow}, and M. Bownik and D. Speegle \cite{BS},  it seems likely that the construction of projection multi-resolution analyses can be generalized to a wider class of dilation matrices of arbitrary dimension. 

%%%%%%%%%%%%%%%%%%%%%%%%%%%%%%%%%%%%%%%%%%%%%%%%%%%%%%%%%%%%%%%%%%%%%%%%%%%%%%%%

%%%%%%%%%%%%%%%%%%%%%%%%%%%
\end{document}